\documentclass{article}
\usepackage{amsmath,amssymb,amsthm,epsfig}
\usepackage[all]{xy}

\usepackage{fullpage}

\theoremstyle{plain}
\newtheorem{proposition}{Proposition}[section]

\theoremstyle{definition}
\newtheorem{definition}{Definition}[section]

\newcommand\A{{\mathcal{A}}}
\newcommand\treek{{X_k^*}}
\newcommand\Ht{{\mathcal H^{(3)}}}
\newcommand\llim{{\underset{n \to \infty}{\lim}}}

\begin{document}

\begin{center}

\vskip 1cm

{\LARGE\bf Tree morphisms, transducers, and \\ \vskip .1in
 integer sequences}

\vskip 1cm \large
Zoran {\v S}uni\'c \footnote{Partially supported by NSF grant DMS-600975}\\
Department of Mathematics\\ Texas A\&M University\\ College Station, TX 77843-3368, USA\\
{\tt sunic@math.tamu.edu}\\
\end{center}

\vskip .2 in

\begin{abstract}
The notion of transducer integer sequences is considered through a
series of examples. By definition, transducer integer sequences are
integer sequences produced, under a suitable interpretation, by
finite automata encoding tree morphisms (length and prefix
preserving transformations of words). Transducer integer sequences
are related to the notion of self-similar groups and semigroups, as
well as to the notion of automatic sequences.
\end{abstract}



\section{Introduction}

It is known from the work of Allouche, B\'etr\'ema, and Shallit
(see~\cite{allouche-b-s:hanoi,allouche-s:as}) that a square free
sequence on 6 letters can be obtained by encoding the optimal
solution to the standard Hanoi Towers Problem on 3 pegs by an
automaton. Roughly speaking, given an input word which is the binary
representation if the number $i$, the automaton ends in one of the 6
states. These states represent the six possible moves between the
three pegs and if the automaton ends in state $q_{xy}$, this means
that the one needs to move the top disk from peg $x$ to peg $y$ in
step $i$ of the optimal solution. The obtained sequence over the
6-letter alphabet $\{\ q_{xy} \mid 0 \leq x,y \leq 2, \ x \neq y \
\}$ is an example of an automatic sequence.

We choose to work with a slightly different type of automata, which
under a suitable interpretation, produce integer sequences in the
output. The difference with the above model, again roughly speaking,
is that not only the final state matters, but the output depends on
every transition step taken during the computation and both the
input and the output words are interpreted as encodings of integers.
The integer sequences that can be obtained this way are called
transducer integer sequences. We provide some examples that
illustrate the notion of a transducer integer sequence. All provided
examples are related to the Hanoi Towers Problem on 3 pegs.

In recent years, a very fruitful line of research in group theory
has led to the notion of a self-similar
group~\cite{nekrashevych:b-selfsimilar} (also known as automata
groups~\cite{grigorchuk-n-s:automata} or state closed
groups~\cite{sidki:circuit}). Many challenging problems have been
solved by using finite automata to encode groups of tree
automorphisms with interesting properties, leading to solutions to
outstanding problems. To name just a few, such examples are the
first Grigorchuk group~\cite{grigorchuk:burnside}, solving the
problem of Milnor on existence of groups of intermediate growth and
the Day-von Neumann problem on existence of amenable but not
elementary amenable groups, Basilica
group~\cite{grigorchuk-z:basilica1,bartholdi-v:basilica}, providing
an example of amenable but not subexponentially amenable group,
Wilson groups~\cite{wilson:nonuniform}, solving the problem of
Gromov on existence of groups of non-uniform exponential growth, the
realization of the lamplighter group $L_2$ by an
automaton~\cite{grigorchuk-z:l2}, leading to the solution of the
Strong Atuyah Conjecture on $L^2$-Betti
numbers~\cite{grigorchuk-al:atiyah}, and the recent solution to
Hubbard's Twisted Rabbit Problem in holomorphic
dynamics~\cite{bartholdi-n:rabbit}. The geometric language and
insight coming from the interpretation of the action of the automata
as tree automorphisms greatly simplifies the presentation and helps
in the understanding of the underlying phenomena, such as
self-similarity, contraction, branching, etc
(see~\cite{grigorchuk-n-s:automata,bartholdi-g-s:branch,bartholdi-g-n:fractal,nekrashevych:b-selfsimilar}
for definitions, examples, and details).

The language of finite automata has been proved extremely suitable
in working with self-similarity phenomena. Indeed, in the case of
automatic sequences, it is known from Cobham
Theorem~\cite{cobham:tags} that such sequences are precisely those
that are obtained as images under codings of fixed points of uniform
endomorphisms (limits of iterations of endomorphisms). The
contracting self-similar groups have been related by Nekrashevych to
finite partial self-coverings of
orbispaces~\cite{nekrashevych:b-selfsimilar}.

In the current article we use automata in the sense of transducers.
As such, they generate self-similar groups (or semigroups) of tree
automorphisms (or endomorphisms). In the same time, the output words
are interpreted as encodings of integers, thus bringing the topic
closer to the topic of automatic sequences. Thus, it is not
surprising that the concrete examples of transducer integer
sequences that are exhibited here all gave high level of
self-similarity and can be defined as limits of certain iterations
of sequences.


\section{Tree morphisms and finite transducers}

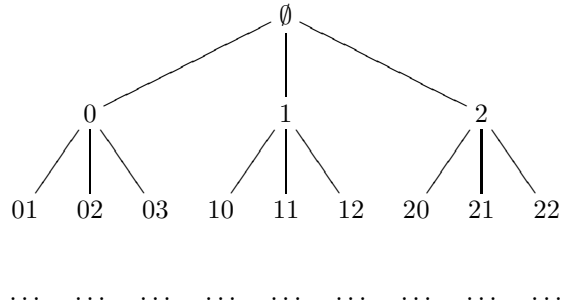
\begin{figure}[!ht]
\begin{center}
\[
 \xymatrix@C=7pt{
  &
  &
  &
  &
  \emptyset \ar@{-}[llld] \ar@{-}[d] \ar@{-}[rrrd]&
  \\
  &
  0 \ar@{-}[ld] \ar@{-}[d] \ar@{-}[rd]&
  &
  &
  1 \ar@{-}[ld] \ar@{-}[d] \ar@{-}[rd]&
  &
  &
  2 \ar@{-}[ld] \ar@{-}[d] \ar@{-}[rd]&
  &
  \\
  01&
  02&
  03&
  10&
  11&
  12&
  20&
  21&
  22
  \\
  \dots&
  \dots&
  \dots&
  \dots&
  \dots&
  \dots&
  \dots&
  \dots&
  \dots&
}
\]
\caption{Ternary rooted tree}\label{3tree}
\end{center}
\end{figure}

For $k \geq 2$ denote $X_k=\{0,1,\dots,k-1\}$. The free monoid
$\treek$ has the structure of a $k$-ary \emph{rooted tree} $\treek$
in which the \emph{empty word} $\emptyset$ is the \emph{root}, the
words of length $n$ constitute the \emph{level $n$} and each vertex
$v$ has $k$ \emph{children}, namely $vx$, for $x$ a letter in $X_k$
(see Figure~\ref{3tree} for the ternary tree).
The tree structure imposes order on $X_k^*$, which is the well known
\emph{prefix order}. Namely, we say that $u \leq v$ if $u$ is a
vertex on the unique geodesic from $v$ to $\emptyset$ in $\treek$,
which is equivalent to saying that $u$ is a \emph{prefix} of $v$. A
map $\mu: X_{k_1}^* \to X_{k_2}^*$ is a \emph{tree morphism} if it
preserves the word length and the prefix relation, i.e
\[ |\mu(u)| = |u| \qquad\text{and}\qquad \mu(u) \leq \mu(uw), \]
for all words $u$ and $w$ over $X_{k_1}$. In case $k_1=k_2$,
morphisms are called \emph{endomorphisms} and bijective
endomorphisms are called \emph{automorphisms}.

Every tree morphism $\mu:X_{k_1}^* \to X_{k_2}^*$ can be decomposed
as
\[ \mu = \pi_\mu (\mu_0,\dots,\mu_{k_1-1}) \]
where $\pi_\mu:X_{k_1} \to X_{k_2}$ is a map called the \emph{root
transformation} of $\mu$ and $\mu_x: X_{k_1}^* \to X_{k_2}^*$, $x$
in $X_{k_1}$, are tree morphisms called the \emph{sections} of
$\mu$. The root permutation and the sections of $\mu$ are uniquely
determined by the recursive relation
\[ \mu(xw) = \pi_\mu(x) \mu_x(w),\]
which holds for every letter $x$ and word $w$ over $X_{k_1}$. Thus
the sections describe the action of $\mu$ on the $k_1$ subtrees
hanging below the root in $X_{k_1}^*$ and the root transformation
$\pi_\mu$ describes the action of $\mu$ at the root.

The tree morphisms act on the left and the composition is performed
from right to left, yielding the formula
\begin{equation}\label{composition}
 \mu \nu = \pi_\mu (\mu_0,\dots,\mu_{k_1-1}) \pi_\nu (\nu_0,\dots,\nu_{k_1-1}) =
           \pi_\mu \pi_\nu (\mu_{\pi_\nu(0)}\nu_0, \dots,
          \mu_{\pi_\nu(k_1-1)}\nu_{k_1-1}).
\end{equation}

A quite efficient way of defining tree morphisms is by using finite
transducers. A \emph{finite $k_1$ to $k_2$ transducer} is a 5-tuple
$\A=(Q,X_{k_1},X_{k_2},\tau,\pi)$, where $Q$ is a finite set of
\emph{states}, $X_{k_1}$ and $X_{k_2}$ are the \emph{input and
output alphabets}, $\tau:Q \times X_{k_1} \to Q$ is a map called the
\emph{transition map} od $\A$, and $\pi:Q \times X_{k_1} \to
X_{k_2}$ is a map called the \emph{output map} of $\A$. Every state
$q$ of the finite transducer $\A$ defines a tree morphism, also
denoted $q$ by setting $q_x=\tau(q,x)$, for $x \in X_{k_1}$, and
$\pi_q: X_{k_1} \to X_{k_2}$ to be the restriction of $\pi$ defined
by $\pi_q(x) = \pi(q,x)$. Thus, for each state $q$ of $\A$ we have
\begin{equation}\label{qxw}
 q(\emptyset) = \emptyset \qquad\text{and}\qquad
 q(xw) = \pi_q(x)q_x(w),
\end{equation}
for $x$ a letter in $X_k$ and $w$ a word over $X_k$. When started at
state $q$, the transducer reads the first input letter $x$, produces
the first letter of the output according to the transformation
$\pi_q$ and changes its state to $q_x$. The state $q_x$ then handles
the rest of the input and output. The states of a $k$-ary transducer
(transducer in which $k_1=k_2=k$) define $k$-ary tree endomorphisms.
An invertible  $k$-ary transducer is a transducer in which
$k_1=k_2=k$ and the transformation $\pi_q$ is a permutation of
$X_k$, for each state $q$ in $Q$. The states of an invertible
$k$-ary transducer define $k$-ary tree automorphisms. When $k_1 \leq
k_2$ and, for each state $q$, the vertex transformation $\pi_q$ is
injective then every state of the transducer $\A$ is an embedding of
the $k_1$-ary tree into the $k_2$-ary tree. We call such a
transducer an \emph{injective transducer}.

The boundary $\partial \treek$ of the $k$-ary tree $\treek$ consists
of all infinite (to the right) words over $X_k$. The boundary has a
structure of an ultrametric space homeomorphic to a Cantor set. The
recursive definition \eqref{qxw} applies to both finite and infinite
words $w$. The action of a state $q$ of a $k$-ary transducer on the
boundary $\partial \treek$ is by continuous maps, while the action
of an invertible $k$-ary transducer is by isometries. More on these
aspects of actions on rooted trees can be found
in~\cite{grigorchuk-n-s:automata}.

\begin{figure}[!ht]
\begin{center}
\begin{tabular}{cc}
 \SelectTips{cm}{}
 \xymatrix@C=45pt@R=45pt{
  &
  *++[o][F-]{{\scriptstyle(01)}}\ar@/_/[d]_{0}\ar@/^/[d]^{1}\ar@(ul,l)_{2}\ar@{}[ddr]^<<<{a_{01}} &
 \\
 & *++[o][F-]{{\scriptstyle()}}\ar@(ul,l)_{0,1,2} \ar@{}[r]^<<<{id}&
 \\
 *++[o][F-]{{\scriptstyle(02)}}\ar@/_/[ur]_{2}\ar@/^/[ur]^{0}\ar@(ul,l)_{1}\ar@{}[rr]_<<<{a_{02}} &
 \A_H &
 *++[o][F-]{{\scriptstyle(12)}}\ar@/_/[ul]_{1}\ar@/^/[ul]^{2}\ar@(ur,r)^{0}\ar@{}[ll]^<<<{a_{12}}
 }
 &
 \SelectTips{cm}{}
  \xymatrix@C=35pt{
  \\ \\
  &
  *++[o][F-]{}
  \ar@/^1pc/[rr]^{2/1} \ar@(u,ul)_{0/0} \ar@(d,dl)^{1/1} \ar@{}[l]|<<<<{\textstyle\alpha}
  &&
  *++[o][F-]{}
  \ar@/^1pc/[ll]^{0/1} \ar@(u,ur)^{1/1} \ar@(d,dr)_{2/0} \ar@{}[r]|<<<<{\textstyle\beta}
  &
  \\
  && \A_L
 }
\end{tabular}
\caption{An invertible ternary transducer $\A_H$ and a 3 to 2
transducer $\A_L$} \label{ah-al}
\end{center}
\end{figure}
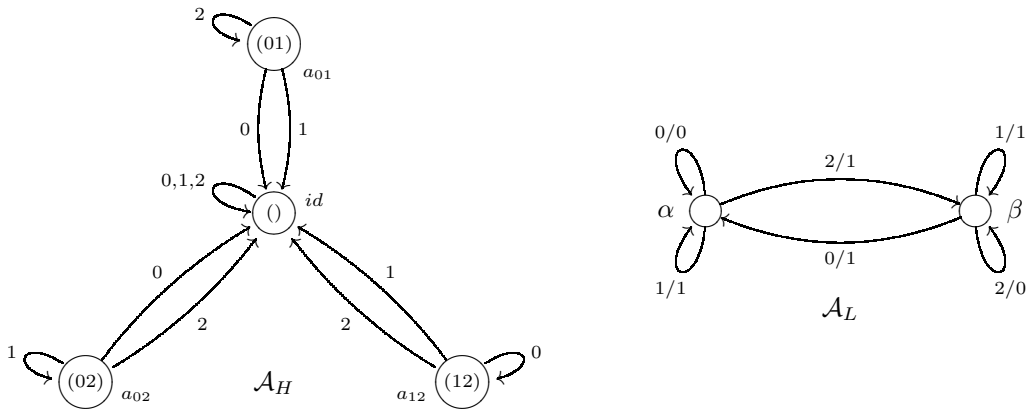

There are two common ways to represent finite $k_1$ to $k_2$
transducers by labeled directed graphs such as the ones in
Figure~\ref{ah-al}. The graph on the left represents an invertible
ternary transducer, denoted $\A_H$. The vertices are the states,
each state $q$ is labeled by its corresponding transformation (in
this case permutation) $\pi_q$, and the edges labeled by the letters
from $X_3$ define the transition function $\tau$ (for every $q$ in
$Q$ and $x$ in $X_3$ there exists an edge from $q$ to
$q_x=\tau(q,x)$ labeled by $x$). The graph on the right represents a
3 to 2 transducer. The vertices are the states and for each pair
$(q,x)$ in $Q \times X_3$ an edge labeled by $x \ | \ \pi_q(x) $
connects $q$ to $q_x$. One can easily switch back and forth between
the two formats. We refer to the second form (the one in which the
output is indicated on the edges) as the \emph{Moore diagram} of the
automaton.

For $0\leq i < j \leq 2$, the ternary tree automorphisms $a_{ij}$
from the automaton $\A_h$ are defined recursively by
\[
\begin{cases}
 a_{ij}(\emptyset) = \emptyset, \\
 a_{ij}(iw) = jw, \\
 a_{ij}(jw) = iw, \\
 a_{ij}(xw) = xa_{ij}(w), & x \not\in \{i,j\},
\end{cases}
\]
for a word $w$ over $X_3$. In simple terms, the only effect the
transformation $a_{ij}$ has on a word $w$ over $X_3$ is that it
changes the very first appearance of either of the symbols $i$ or
$j$ in $w$ to the other symbol, if such an appearance exists. To
simplify the notation, we write
\[ a = a_{01}, \qquad b = a_{02}, \qquad\text{and}\qquad c = a_{12}. \]
The state labeled by $id$ does not change any input word and
represents the identity automorphism of the ternary tree. It is
clear that $a$, $b$ and $c$ are self-invertible transformations of
$X_3^*$, i.e~$a^2=b^2=c^2=id$.

The 3 to 2 tree morphisms defined by the transducer $A_L$ are
defined recursively by
\begin{alignat*}{4}
 \alpha(\emptyset) &= \emptyset, \qquad & \alpha(0w) &= 0\alpha(w),\qquad &
 \alpha(1w) &= 1\alpha(w), \qquad & \alpha(2w) &= 1\beta(w), \\
 \beta(\emptyset) &= \emptyset, \qquad & \beta(0w) &= 1\alpha(w),\qquad &
 \beta(1w) &= 1\beta(w), \qquad & \beta(2w) &= 0\beta(w).
\end{alignat*}

\begin{definition}
The semigroup (group) of $k$-ary tree endomorphisms (automorphisms)
generated by all the states of an (invertible) $k$-ary transducer
$\A$ is called the semigroup (group) of $\A$ and is denoted by
$S(\A)$ ($G(\A)$).
\end{definition}

The group $G(\A_H)$ is introduced in~\cite{grigorchuk-s:hanoi-cr},
where it is called Hanoi Towers group on 3 pegs and denoted $\Ht$
(in fact, one Hanoi Towers group $\mathcal{H}^{(k)}$ is introduced
for each number of pegs $k \geq 3$). The name is derived from the
fact that the group $\Ht$ models the well known Hanoi Towers Problem
on 3 pegs.

To recall, the Hanoi Towers Problem on 3 pegs and $n$ disks is the
following. In a valid $n$ disk configuration, disks of different
size, labeled by $1,2,\dots,n$ according to their size, are placed
on three pegs, labeled 0,1 and 2, in such a way that no disk is
placed on top of a smaller disk. In a single move the top disk from
one peg can be moved and placed on top of another peg, as long as
the newly obtained configuration is still valid. Initially all $n$
disks are placed on peg 0 and the problem asks for an optimal
algorithm that moves all disks to another peg.

\begin{figure}[!ht]
\begin{center}
\epsfig{file=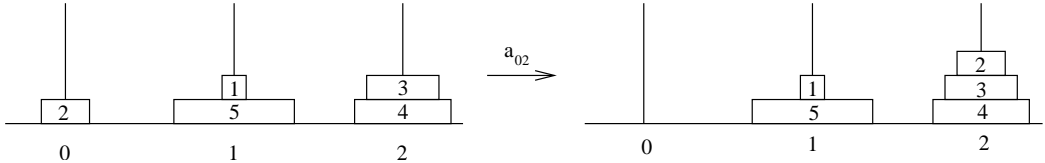,height=60pt} \caption{A move between peg 0
and peg 2}\label{move}
\end{center}
\end{figure}

Each valid configuration of $n$ disks can be encoded by a word of
length $n$ over $X_3$. The word $x_1 \dots x_n$ represents the
unique valid configuration in which disk $i$ is placed on peg $x_i$.
The ternary tree automorphism $a_{ij}$ then represents a move
between peg $i$ and peg $j$ (in either direction). For example the
move between peg 0 and peg 2 illustrated in Figure~\ref{move} is
encoded as $a_{02}(10221) = 12221$.

The action of $\Ht$ on the ternary tree is spherically transitive,
meaning that it is transitive on the levels of the tree. This is
equivalent to the statement that any valid configuration on $n$
disks can be obtained from any other valid configuration on $n$
disks by legal moves.

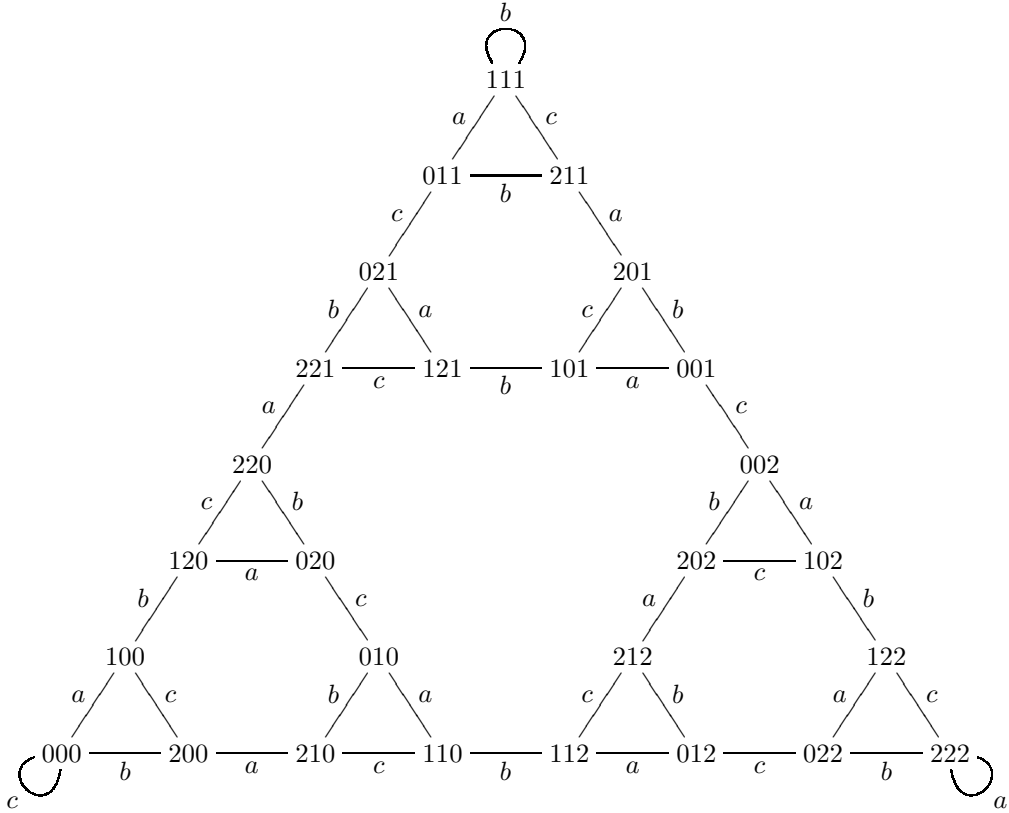
\begin{figure}[!ht]
\begin{center}
 \[
 \SelectTips{cm}{}
 \xymatrix@C=3pt{
 & & & & &  & & 111 \ar@{-}@(ur,ul)_{\textstyle b}
 \\
 & & & & &  & 011 \ar@{-}[rr]_{\textstyle b} \ar@{-}[ur]^{\textstyle a} & & 211 \ar@{-}[ul]_{\textstyle c}
 \\
 & & & & & 021 \ar@{-}[ur]^{\textstyle c} & & & & 201 \ar@{-}[ul]_{\textstyle a}
 \\
 &  & & & 221 \ar@{-}[rr]_{\textstyle c} \ar@{-}[ur]^{\textstyle b}
      & & 121 \ar@{-}[rr]_{\textstyle b} \ar@{-}[ul]_{\textstyle a}
      & & 101 \ar@{-}[rr]_{\textstyle a} \ar@{-}[ur]^{\textstyle c}
      & & 001 \ar@{-}[ul]_{\textstyle b}
 \\
 & & & 220 \ar@{-}[ur]^{\textstyle a} & & & & & & & & 002 \ar@{-}[ul]_{\textstyle c}
  \\
 & & 120 \ar@{-}[rr]_{\textstyle a} \ar@{-}[ur]^{\textstyle c} & & 020 \ar@{-}[ul]_{\textstyle b} & & & &
 & & 202 \ar@{-}[rr]_{\textstyle c} \ar@{-}[ur]^{\textstyle b} & & 102 \ar@{-}[ul]_{\textstyle a}
 \\
 & 100 \ar@{-}[ld]_{\textstyle a} \ar@{-}[rd]^{\textstyle c} \ar@{-}[ur]^{\textstyle b} &&&&
   010 \ar@{-}[ld]_{\textstyle b} \ar@{-}[rd]^{\textstyle a} \ar@{-}[ul]_{\textstyle c} &&&&
   212 \ar@{-}[ld]_{\textstyle c} \ar@{-}[rd]^{\textstyle b} \ar@{-}[ur]^{\textstyle a} &&&&
   122 \ar@{-}[ld]_{\textstyle a} \ar@{-}[rd]^{\textstyle c} \ar@{-}[ul]_{\textstyle b}
 \\
 000 \ar@{-}[rr]_{\textstyle b} \ar@{-}@(l,d)_{\textstyle c}  &&
 200 \ar@{-}[rr]_{\textstyle a} && 210 \ar@{-}[rr]_{\textstyle c} &&
 110 \ar@{-}[rr]_{\textstyle b} && 112 \ar@{-}[rr]_{\textstyle a} &&
 012 \ar@{-}[rr]_{\textstyle c} && 022 \ar@{-}[rr]_{\textstyle b} &&
 222 \ar@{-}@(r,d)^{\textstyle a}
 }
 \]
\caption{The Schreier graph of $\Ht$ at level 3}\label{sierpinski3}
\end{center}
\end{figure}

Consider the stabilizer of the vertex $0^n$ in $\Ht$, denoted $P_n$.
The group $\Ht$ acts on the set $\Ht/P_n$ of left cosets of $P_n$.
The action is described by the corresponding Schreier graph
$\Gamma_n = \Gamma_n(\Ht,P_n,S)$ of $P_n$ with respect to the
generating set $S=\{a,b,c\}$. The vertices are the cosets of $P_n$
and there is an edge connecting $hP_n$ to $shP_n$ for every coset
$hP_n$ and generator $s$ in $S$. Since $h' \in hP_n$ if and only if
$h'(0^n) = h(0^n)$ the vertices of the Schreier graph $\Gamma_n$ can
be encoded by the vertices of the $n$-h level of the ternary tree
(the coset $hP_n$ os labeled by $h(0^n)$) and two vertices are
connected if and only if one is the image of the other under $s$,
for some generator $s$ in $S$. The Schreier graph $\Gamma_3$
corresponding to level 3 of the ternary tree is given in
Figure~\ref{sierpinski3}.
Since all generators have order 2, no directions are indicated on
the edges.

The sequence of graphs $\{\Gamma_n\}$ converges to an infinite graph
$\Gamma$ in the space of pointed graphs based at $0^n$
(see~\cite{grigorchuk-z:cortona} for definitions of this space),
which is the Schreier graph $\Gamma=\Gamma(\Ht,P,S)$, where
$P=\cap_{n=0}^\infty P_n$ is the stabilizer of the infinite ray
$0^\infty=000\dots$ on the boundary of the ternary tree. One can
think of the limiting graph both as the Schreier graphs of the
action of $\Ht$ on the orbit of the infinite ray $0^\infty$ in
$\partial X_3^*$ or as the model of Hanoi Towers Problem
representing all valid configurations that can be reached from the
configuration in which (countably) infinitely many disks are placed
on peg $0$ (this configuration corresponds to the infinite word
$0^\infty$)

Graphs similar to $\Gamma_n$, modeling the Hanoi Towers problem are
well known in the literature, but there is a subtle difference.
Namely, the difference with the corresponding graphs
in~\cite{hinz:ens} modeling the Hanoi Towers Problem is that the
edges in $\Gamma_3$ are labeled (by the corresponding tree
automorphisms) and our graphs have loops at the corners
(corresponding to situations in which all disks are on one peg and
the generator corresponding to a move between the other two pegs
does not change anything), which turn them into 3-regular graphs.
Finite dimensional permutational representations of $\Ht$ based on
the action on the levels of the ternary tree were used
in~\cite{grigorchuk-s:hanoi-cr} to calculate the (Markov) spectrum
of the graphs $\Gamma_n$ as well as the limiting infinite graph
$\Gamma$. Among interesting properties of $\Ht$ we mention that it
is an amenable (but not subexponentially amenable), regular branch
group over its commutator, it is not just infinite and its closure
in the pro-finite group of ternary tree automorphisms is finitely
constrained. Moreover, $\Ht$ is (up to conjugation) the iterated
monodromy group of the finite rational map $z \mapsto z^2 -
\frac{16}{27z}$, whose Julia set is the Sierpi{\'n}ski gasket. This
explains the fact that the sequence of Schreier graphs
$\{\Gamma_{0^n}\}$ approximates the Sierpi{\'n}ski gasket. For more
information on properties of $\Ht$ we refer the interested reader
to~\cite{grigorchuk-s:hanoi-cr,grigorchuk-s:standrews,grigorchuk-n-s:oberwolfach1,grigorchuk-n-s:oberwolfach2}.

When $k_1 \geq k_2$ every $k_1$ to $k_2$ tree morphism can also be
thought of as a $k_1$-ary tree morphism, since the $k_2$-ary tree
canonically embeds in the $k_1$-ary tree in obvious way. We
calculate the semigroup $S(\A_L)$ by thinking of the transducer
$\A_L$ as being a ternary transducer.

\begin{proposition}
The self-similar semigroup $S_L=S(\A_L)$ is given by the semigroup
presentation
\[
 S_L = \langle \alpha, \beta \mid \alpha^2=\alpha, \ \alpha\beta = \beta \ \rangle.
\]
In other words, $S(\A_L)$ is the free cyclic semigroup generated by
$\beta$ extended by a left identity element $\alpha$.
\end{proposition}

\begin{proof}
Since $\alpha$ acts trivially on the binary words (words over
$\{0,1\}$) and the image of every ternary word under the elements of
$S_L$ is a binary word, we have $\alpha \sigma = \sigma$, for every
element $\sigma$ of $S_L$.

Denote by $\pi_1$ the transformation $X_3 \to X_3$ given by
$\pi_1(x)=1$, for $x$ in $X_3$. Note that $\pi_\beta\pi_\sigma =
\pi_1$, for all elements $\sigma$ of $S_L$. Calculations using
\eqref{composition} yield
\begin{align*}
 \beta   &= \pi_\beta ( \alpha, \beta, \beta), \\
 \beta^2 &= \pi_1     ( \beta\alpha, \beta\beta, \alpha\beta) =
            \pi_1( \beta\alpha, \beta^2, \beta), \\
\beta^3 &= \pi_1     ( \beta\beta\alpha, \beta\beta^2, \beta\beta) =
             \pi_1( \beta^2\alpha, \beta^3, \beta^2) \\
\dots \\
\beta^n &= \pi_1 (\beta^{n-1}\alpha, \beta^n,\beta^{n-1}), \qquad
\text{for }n \geq 2.
\end{align*}
Since $\pi_\beta \neq \pi_1$ the tree morphism $\beta$ is different
than any tree morphism $\beta^n$, for $n \geq 2$. On the other hand
if $\beta^{n_1} = \beta^{n_2}$, for some $n_1,n_2\geq 2$ then their
sections at coordinate 2 must be equal, which forces $\beta^{n_1-1}
= \beta^{n_2-1}$. Finite descent then shows that all positive powers
of $\beta$ are distinct.

Further,
\begin{align*}
 \alpha   &= \pi_\alpha ( \alpha, \alpha, \beta), \\
 \beta^m\alpha &= \pi_1 (\beta^{m-1}\alpha, \beta^m\alpha,\beta^{m+1}), \qquad
\text{for }m \geq 1.
\end{align*}
The powers of $\beta$ in coordinate 2 imply that all the elements
$\beta^m\alpha$ are distinct for distinct values of $m$.

Finally, assuming $\beta^m\alpha = \beta^n$, for some $m$ and $n$,
forces $2 \leq n = m+2$, by comparing the sections at coordinate 2.
However, $\alpha \neq \beta^2$ since they have different root
transformation. For $m \geq 1$ the equality $\beta^m\alpha =
\beta^{m+2}$ implies $\beta^{m-1}\alpha = \beta^{m+1}$, by comparing
the sections at coordinate 1. Finite descent then finishes the
proof.
\end{proof}


\section{Transducer integer sequences}

We first recall the well established notion of automatic sequence.
The definition that follows is one of the equivalent definitions
that can be found in~\cite{allouche-s:as}.

A $k$-ary \emph{finite automaton with final state output} ($k \geq
2$) is a 6-tuple $\A=(Q,X_k,Y,s,\tau,\pi)$, where $Q$ is a finite
set, called set of \emph{states}, $X_k=\{0,\dots,k-1\}$ is the
\emph{input alphabet}, $Y$ is a finite set called the \emph{output
alphabet}, $s$ is an element in $Q$ called the \emph{initial state},
$\tau:Q \times X \to Q$ is a map called \emph{transition map} and
$\pi:Q \to Y$ is a map called \emph{final state output map}. Such an
automaton defines an infinite sequence $y_0,y_1,y_2,\dots$ over the
output alphabet $Y$, called the \emph{final state output sequence}
of $A$, as follows. For a natural number $i\geq 0$ let
$[i]_k=i_0\dots i_m$ be any base $k$ representation of $i$ with $i =
\sum_{j=1}^m i_j k^j$ (thus the least significant digit is written
first). The term $y_i$ in the final state output sequence is defined
as the image $\pi(q)$ of the state $q$ the automaton reaches as it
reads the input word $[i]_k$ starting from the initial state $s$
(this output must be independent of the chosen representation of
$i$). Thus
\[ y_i = \pi(\tau(s,[i]_k)), \]
where $\tau:Q \times X^* \to Q$ is the recursive extension of $\tau$
on $Q \times X^*$ defined by $\tau(q,\emptyset)=q$ and $\tau(q,xw) =
\tau(\tau(q,x),w)$, for $q$ a state in $Q$, $x$ a letter in $X_k$
and and $w$ a word over $X_k$.

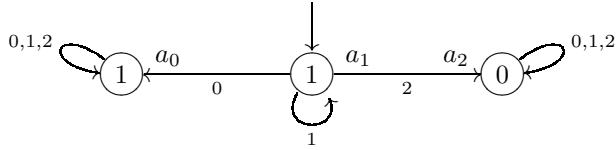
\begin{figure}[!ht]
\begin{center}
 \[
  \SelectTips{cm}{}
 \xymatrix@C=55pt{
  *++[o][F-]{1} \ar@(ul,l)_{0,1,2} \ar@{}[r]^<<<<{\textstyle a_0}
  &
  *++[o][F-]{1}\ar@{->}[l]^{0}\ar@{->}[r]_{2}\ar@(dl,dr)_{1}\ar@(u,u)\ar@{}[r]^<<<<{\textstyle a_1}
  &
  *++[o][F-]{0} \ar@(ur,r)^{0,1,2}\ar@{}[l]_<<<<{\textstyle a_2}
 }
 \]
\caption{A ternary automaton with final state output
$A_{0-2}$}\label{a02}
\end{center}
\end{figure}

Automata with final state output can be represented by labeled
directed graphs similar to the ones representing transducers. The
only significant difference is that each state $q$ is labeled by the
corresponding output letter $\pi(q)$ and the initial state is
indicated by an incoming arrow. As an example, consider the ternary
automaton $\A_{0-2}$ in Figure~\ref{a02}.

\begin{definition}
A $k$-ary automatic sequence is an infinite sequence that can be
obtained as the final state output sequence of some $k$-ary finite
automaton.
\end{definition}

By Cobham Theorem~\cite{cobham:tags} a sequence over a finite
alphabet is a $k$-ary automatic sequence if and only if it is an
image under a coding of a fixed point of a $k$-uniform endomorphism.

Given a free monoid $X^*$ over a finite alphabet $X$, an
endomorphism $\alpha:X^* \to X^*$ can be uniquely defined by
specifying the images of the letters in $X$ under $\alpha$. Let
there exists a letter $x$ in $X$ such that $\alpha(x) = xw$, where
$w$ is non-empty word, and let $\alpha(x)\neq \emptyset$, for all
letters $x$ in $X$. Then, for all $n \geq 0$ the $n$-th iterate
$\alpha^n(x)$ is a proper prefix of the $(n+1)$-st iterate
$\alpha^{n+1}(x)=\alpha(\alpha^n(x))$ and the limit $\llim
\alpha^n(x)$ is a well defined infinite sequence over $X$. In the
particular case when the length of all the words $\alpha(x)$, $x \in
X$, is equal to $k$, the morphism $\alpha$ is called a $k$-uniform
endomorphism.

As an example, let $X=\{1,-1\}$ and denote by $w_\alpha$ the
infinite binary sequence
\[ w_\alpha = \llim \alpha^n(1) = 11\text{-}1 \ 11\text{-}1 \ 1\text{-}1\text{-}1 \ 11\text{-}1 \ 11\text{-}1 \ 1\text{-}1\text{-}1 \ 11\text{-}1 \ 1\text{-}1\text{-}1 \ 1\text{-}1\text{-}1 \dots~ \]
obtained by iterations, starting from 1, of the endomorphism
$\alpha:X^* \to X^*$ given by (compare to the sequence A080846)
\[ 1 \mapsto 11\text{-}1 \qquad \text{-}1 \mapsto 1\text{-}1\text{-}1. \]

A finite or infinite word $w$ over an alphabet $X$ is cube free if
it does not contain a subword of the form $uuu$, where $u$ is a
nontrivial finite word over $X$.

\begin{proposition}
The infinite binary sequence $w_\alpha$ is cube-free.
\end{proposition}
\begin{proof}
By the criterion of Richomme and
Wlazinski~\cite{richomme-w:cube-free}, an easy way to verify that
$w_\alpha$ is cube free is to observe that
$\alpha(11\text{-}1\text{-}11\text{-}11\text{-}1\text{-}11\text{-}1\text{-}111\text{-}111\text{-}11\text{-}111\text{-}1\text{-}1)$
is cube free.
\end{proof}

We offer two additional descriptions of $w_\alpha$.

Define a sequence of words $w_{[n]}$ of length $3^n$ by
\begin{align*}
 w_{[0]}     &= 1, \\
 w_{[n+1]} &= w_{[n]} w_{[n]} w_{[n]}',
\end{align*}
where $w_{[n]}'$ is obtained from $w_{[n]}$ by changing the middle
symbol in $w_{[n]}$ from 1 to -1.

\begin{proposition}
The limit $\llim w_{[n]}$ is well defined and is equal to
$w_\alpha$.
\end{proposition}

For an integer $i \geq 0$, let $(i)_k=i_0i_1 \dots $ be the sequence
of digits in base $k$ representation of $i$, where $i =
\sum_{j=0}^\infty i_jk^j$ (the sequence ends in infinitely many
0's).

Call a natural number $i$ a 2-before-0 number if the least
significant digit in the ternary representation $(i)_3$ of $i$ that
is different from 1 is 2. Otherwise the number is called a
0-before-2 number. Define an infinite binary sequence
$x_0,x_1,x_2,\dots,$ by
\[ x_i = \begin{cases}
      1, & \text{if }i \text{ is a 0-before-2 number}\\
     -1, & \text{if }i \text{ is a 2-before-0 number}
          \end{cases}.
\]

\begin{proposition}
The infinite binary sequence $x_0,x_1,x_2,\dots,$ is equal to
$w_\alpha$.
\end{proposition}

\begin{proposition}
The infinite binary sequence $w_\alpha$ is a ternary automatic
sequence. It can be obtained as the final state output sequence of
the automaton $A_{0-2}$.
\end{proposition}
\begin{proof}
The only time the automaton $A_{0-2}$ produces -1 in the output is
if it reaches the state $a_2$, which only happens if $i$ is a
2-before-0 number.
\end{proof}

We define now the notion of transducer integer sequence.

\begin{definition}
A $k_1$ to $k_2$ \emph{transducer integer sequence} is a sequence of
integers $\{z_i\}_{i=0}^\infty$ such that there exists a $k_1$ to
$k_2$ transducer $\A$ and a state $q$ in $\A$ such that, for every
$i \geq 0$, the output word $q((i)_{k_1})$ is the base $k_2$
representation of $z_i$.
\end{definition}

It is implicit in the above definition that the state $q$ of $\A$
maps the confinality class of $0^\infty$ in $\partial X_{k_1}^*$ to
the confinality class of $0^\infty$ in $\partial X_{k_2}^*$ (the
confinality class of $0^\infty$ is just the set of infinite words
ending in $0^\infty$). We keep our attention only to this class
since it is the one describing non-negative integers.

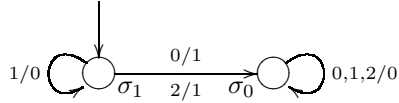
\begin{figure}[!ht]
\begin{center}
\[
 \xymatrix{
 *++[o][F-]{} \ar@{->}[rr]^{0/1}_{2/1} \ar@(ul,dl)_{1/0} \ar@(u,u) \ar@{}[rr]_<<<{\textstyle \sigma_1} &&
 *++[o][F-]{} \ar@(ur,dr)^{0,1,2/0} \ar@{}[ll]^<<<{\textstyle \sigma_0}
 }
\]
\caption{A ternary transducer $\A_T$}\label{at}
\end{center}
\end{figure}

As an easy example, let $\A_T$ be the ternary transducer in
Figure~\ref{at}.
The state labeled by $\sigma_0$ just rewrites all digits to 0.
Clearly
\[ \sigma_1(1^n0w) = \sigma_1(1^n2w) = 0^n10^\infty \]
for any word $w$ in the confinality class of $0^\infty$. Since
$[0^n1]_3 = 3^n$ the obtained integer sequence
$\{a_n\}_{n=0}^\infty$ is (compare to sequence A038500)
\[ 1,3,1,\ 1,9,1,\ 1,3,1,\ 1,3,1,\ 1,27,1,\ 1,3,1,\ 1,3,1,\ 1,9,1,\ 1,3,1,\dots~. \]
By thinking of the powers of 3 as an (infinite) alphabet, this
sequence can be thought of as the fixed point of the iterations
starting from 1 of the 3-uniform endomorphism defined by
\[ x \mapsto 1, \ 3x, 1~. \]
This sequence can also be defined by blocks $a_{[n]}$ of length
$3^n$ as
\[ a_{[0]} = 1 \qquad\qquad a_{[n+1]} = a_{[n]}a_{[n]}'a_{[n]}, \]
where $a_{[n]}'$ is obtained from $a_{[n]}$ by multiplying the
middle term by 3.

A more interesting example is provided by the automaton $\A_L$.

Let $N_2$ be the set of all non-negative integers whose base 3
representation does not use the digit 2 (they are listed in sequence
A005836). Define a sequence $\{\ell_n\}_{n=0}^\infty$, called
\emph{L-sequence},  by
\[ \ell_n = \ell_n^+ + \ell_n^- \]
where $\ell_n^+$ and $\ell_n^-$ are the unique non-negative integers
such that $\ell_n^+,\ell_n^-,\ell_n^+ +\ell_n^- \in N_2$ and
$n=\ell_n^+ - \ell_n^-$ (sequence A060374).

\begin{proposition}
The $L$-sequence is a ternary transducer integer sequence. It is
generated by the transducer $A_L$ with initial state $\alpha$.
\end{proposition}
\begin{proof}
When the current input digit of $n$ is 0, the corresponding digit in
$\ell_n^-$ must be 0. Indeed if it were 1, then the corresponding
digit in $\ell_n^+$ would be $0+1=1$, which would force the
corresponding digit in $\ell_n$ to be 2. Thus the corresponding
digit in $\ell_n^-$ is 0, and so are the digits in $\ell_n^+$ and
$\ell_n$. This corresponds to the first column under $\alpha$ in the
following table.
\[
 \begin{array}{c|ccc|ccc}
          &   & \alpha &&& \beta \\
 \hline
 n        & 0 & 1 & 2 & 0 & 1 &2 \\
 \ell_n^- & 0 & 0 & 1 & 0 & 1 &0  \\
 \ell_n^+ & 0 & 1 & 0 & 1 & 0 &0 \\
 \hline
 \ell     & 0 & 1 & 1*& 1 & 1*&0*\\
 \end{array}
\]
Similarly, if the current input digit in $n$ is 1, the corresponding
digit in $\ell_n^-$ must be 0, in $\ell_n^+$ must be 1 and in
$\ell_n$ must be 1. If the current input digit in $n$ is 2 then the
corresponding digit in $\ell^-$ must be 1. Indeed if it were $0$
then the corresponding digit in $\ell_n^+$ would be $2+0=2$. Thus
the corresponding digit in $\ell_n^-$ is 1, in $\ell_n^+$ is 0 and
in $\ell_n$ is 1. However, in this case there is a carryover for the
next digit (indicated by the $*$ in the table. This is why a second
state $\beta$ is introduced (this state corresponds to the situation
in which there is a carryover in the addition
$n+\ell_n^-=\ell_n^+$). The entries in the right half of the table
(those corresponding to $\beta$) can be treated similarly, by taking
into account the carryover.
\end{proof}

Let $\{p_n\}_{n=0}^\infty$ be the sequence defined by
\[
 p_0 = 0, \qquad\qquad
 p_n= \sum_{i=0}^{n-1} w_ia_i, \quad\text{ for }n \geq 1
\]
where the sequence $\{w_n\}_{n=0}^\infty$ providing the signs is the
cube free sequence generated by the automaton $\A_{0-2}$ and
$\{a_n\}_{n=0}^\infty$ is the transducer sequence generated by
$\A_T$.

\begin{proposition}
The sequence $\{p_n\}$ is equal to the $L$-sequence.
\end{proposition}
\begin{proof}
We have $p_0=0=\ell_0$ and, for $n$ a positive integer and $w$ a
word over $X_3$,
\begin{align*}
  \alpha(0w+1) &= \alpha(1w) = 1\alpha(w) = 0\alpha(w)+1 = \alpha(0w)+1,\\
  \alpha(1^n0w+1) &= \alpha(21^{n-1}0w) = 11^{n-1}1\alpha(w) =
   1^n0\alpha(w) + 3^n = \alpha(1^n0w) + 3^n,\\
  \alpha(1^n2w+1) &= \alpha(21^{n-1}2w) = 11^{n-1}0\beta(w) =
   1^n1\beta(w) - 3^n = \alpha(1^n2w) - 3^n, \\
  \alpha(2^n0w+1) &= \alpha(0^n1w) = 0^n1\alpha(w) =
   10^{n-1}1\alpha(w) - 1 = \alpha(2^n0w) - 1, \\
  \alpha(2^n1w+1) &= \alpha(0^n2w) = 0^n1\beta(w) =
   10^{n-1}1\alpha(w) - 1 = \alpha(2^n1w) - 1. \\
\end{align*}
In each case the change in the value of $\alpha(i)$ is exactly
$w_ia_i$, i.e.~, for all $i$,
\[ \ell_{i+1} = \alpha(i+1) = \alpha(i) + w_ia_i = \ell_i+ w_ia_i \]
and therefore the sequence of partial sums $\{p_n\}$ is exactly the
$L$-sequence.
\end{proof}

The sequence $\{\ell_n\}_{n=0}^\infty$ can also be described as a
fixed point of an endomorphism over the alphabet consisting of the
elements of $N_2$. The iterations start at 0 and the endomorphism is
given by
\[
   0 \mapsto 0, 1 \qquad\qquad
   x \mapsto 3x+1, 3x, 3x+1, \quad \text{ for } x\geq 1.
\]

\section{Relation to Hanoi Towers Problem}

In this section we exhibit a connection between Hanoi Towers
Problem, the automatic cube free sequence $\{w_n\}$ and the
transducer sequence $\{a_n\}$.

Define a matrix $K_n$ of size $3^n \times n$ with entries in $X_3$
by
\[
 K_1 = \begin{bmatrix} 0 \\ 1 \\ 2 \end{bmatrix}, \hspace{3cm}
 K_{n+1} = \begin{bmatrix} K_n & 0_n \\ K_n^R & 1_n \\ K_n & 2_n \end{bmatrix},
\]
where the matrix $K_n^R$ is obtained from the matrix $K_n$ by
flipping $K_n$ along the horizontal axis, and $0_n$, $1_n$ and $2_n$
are column vectors with $3^n$ entries equal to 0, 1 and 2,
respectively. Denote the infinite limit matrix $\underset{n \to
\infty}{\lim} K_n$ by $K$.

For example, the transpose of $K_3$ is given by
\[
 K_3^T =
\begin{bmatrix}
 012 & 210 & 012 & 210 & 012 & 210 & 012 & 210 & 012 \\
 000 & 111 & 222 & 222 & 111 & 000 & 000 & 111 & 222 \\
 000 & 000 & 000 & 111 & 111 & 111 & 222 & 222 & 222
\end{bmatrix}
\]

The limiting matrix $K$ is well defined due to the fact that $K_n$
appears as the upper left corner in $K_{n+1}$. By definition, the
indexing of the rows of $K$ starts with 0 while the indexing of the
columns starts with 1.

A sequence $w_0,\dots,w_{k^n-1}$ of words of length $n$ over $X_k$
is a $k$-ary Gray code of length $n$ if all words of length $n$ over
$X_k$ appear exactly once in the sequence and any two consecutive
words differ in exactly one position.

\begin{proposition}
The $3^n$ rows of the matrix $K_n$ represent a ternary Gray code of
length $n$.
\end{proposition}

By interpreting the rows of $K$ as ternary representations of
integers, we obtain the sequence
\[ 0,1,2,5,4,3,6,7,8,17,16,15,12,13,14,11,10,9,\dots, \]
which is not included in The On-Line Encyclopedia of Integer
Sequences (as of December 2006).

We observe that the successive rows in $K$ are obtained from each
other by applying the ternary tree automorphism $a$ at odd steps and
$c$ at even steps (the automorphisms $a$ and $c$ are defined by
$\A_H$ - the automaton generating the Hanoi Towers group).

\begin{proposition}
For $j \geq 0$, define $t_{2j} = (ca)^j$ and $t_{2j+1} = a(ca)^j$.
Let $k_i$ denote the $i$-th row in the matrix $K$. Then
\[ k_i = t_i(k_0). \]
\end{proposition}
\begin{proof}
The proof is by induction on $n$ in $K_n$. The crucial observation
for the inductive step is that the last row in $K_n$ is
$2^n0^\infty$ and is obtained by applying $c$ in step $3^n-1$. In
the next step applying $a$ to $2^n0^\infty$ produces $2^n1^\infty$.
Alternate applications of $c$ and $a$ ($3^n-1$ total) do not affect
the $1$ in the position $n+1$, but backtrack the word in the first
$n$ entries from $2^n$ back to $0^n$, thus producing $0^n10^\infty$
at step $3^n-1+1+3^n-1=2\cdot 3^n-1$. The last taken step is $a$ so
$c$ takes $0^n10^\infty$ to $0^n20^\infty$ and then alternate
applications of $a$ and $c$ change the first $n$ entries again from
$0^n$ to $2^n$ in $3^n-1$ steps, eventually producing
$2^{n+1}0^\infty$ in $2\cdot 3^n-1 + 1 + 3^n-1 = 3^{n+1}-1$  steps,
alternating between $a$ and $c$.
\end{proof}

It is clear that the rows of $K_i$ constitute the whole confinality
class of $0^\infty$. Thus the subgroup $\langle a, c\rangle$ acts
transitively on this class. Since the order of both $a$ and $c$ is 2
this means that $\langle a,c\rangle$ is the infinite dihedral group
$D_\infty$. The transitivity of the action of $\langle a,c\rangle$
on the confinality class of $0^\infty$ is equivalent to the known
fact that any valid $n$ disk configuration can be obtained from any
other in a restricted version of  Hanoi Towers Problem in which no
disk can move between pegs 0 and 2 (in our terminology, applications
of the automorphism $b$ are not allowed). Figure~\ref{peano} shows
the path taken by $(ca)^{13}$ from $000$ to $222$ in $\Gamma_3$.
\begin{figure}[!ht]
\begin{center}
 \[
 \SelectTips{cm}{}
 \xymatrix@C=3pt{
 & & & & &  & & 111
 \\
 & & & & &  & 011 \ar@{-}[ur]^{\textstyle a} & & 211 \ar@{-}[ul]_{\textstyle c}
 \\
 & & & & & 021 \ar@{-}[ur]^{\textstyle c} & & & & 201 \ar@{-}[ul]_{\textstyle a}
 \\
 &  & & & 221 \ar@{-}[rr]_{\textstyle c} & & 121  \ar@{-}[ul]_{\textstyle a}
      & & 101 \ar@{-}[rr]_{\textstyle a} \ar@{-}[ur]^{\textstyle c} & & 001
  \\
 & & & 220 \ar@{-}[ur]^{\textstyle a} & & & & & & & & 002 \ar@{-}[ul]_{\textstyle c}
  \\
 & & 120 \ar@{-}[rr]_{\textstyle a} \ar@{-}[ur]^{\textstyle c} & & 020  & & & &
 & & 202 \ar@{-}[rr]_{\textstyle c}  & & 102 \ar@{-}[ul]_{\textstyle a}
 \\
 & 100 \ar@{-}[ld]_{\textstyle a} \ar@{-}[rd]^{\textstyle c}  &&&&
   010  \ar@{-}[rd]^{\textstyle a} \ar@{-}[ul]_{\textstyle c} &&&&
   212 \ar@{-}[ld]_{\textstyle c}  \ar@{-}[ur]^{\textstyle a} &&&&
   122 \ar@{-}[ld]_{\textstyle a} \ar@{-}[rd]^{\textstyle c}
 \\
 000  &&  200 \ar@{-}[rr]_{\textstyle a} && 210 \ar@{-}[rr]_{\textstyle c} &&
 110  &&  112 \ar@{-}[rr]_{\textstyle a} && 012 \ar@{-}[rr]_{\textstyle c} &&
 022  &&  222
 }
 \]
\caption{The ternary Gray code path generated by $a$ and $c$ in
$\Ht$ at level 3}\label{peano}
\end{center}
\end{figure}

Order all configurations (words in the confinality class of
$0^\infty$) according to their position in the matrix $K$ (small
configurations correspond to rows with small index). When $b$ is
applied to any configuration $k_i$ the obtained configuration
$b(k_i)$ is either larger or smaller than $k_i$. Based on this
alternative define an infinite sequence $\{d_i\}_{i=0}^\infty$ over
$X=\{1,-1\}$ by
\[ d_i = \begin{cases}
               1, & \text{if } b(k_i) > k_i \\
              -1, & \text{if } b(k_i) < k_i
          \end{cases}.
\]
Call this sequence the $b$-direction sequence. Further, define an
integer sequence $\{b_i\}_{i=0}^\infty$ by $b_i = |i-j|$, where $j$
is the index of the configuration $k_j=b(k_i)$. Call this sequence
the $b$-change index sequence.

\begin{proposition}
The $b$-direction sequence is exactly the cube free automatic
sequence $\{w_n\}$ generated by $\A_{0-2}$ and the $b$-change index
sequence is exactly the transducer integer sequence $\{a_n\}$
generated by $\A_T$.
\end{proposition}
\begin{proof}
The proof is by induction on blocks of size $3^n$. Observe that in
each matrix $K_n$ the configuration which is half way between $0^n$
and $2^n$ is $1^n$. The size $3^n$ blocks of the $b$-change sequence
satisfy a relation of the form $b_{[n+1]} =
b_{[n]}'b_{[n]}''b_{[n]}'''$, where $b_{[n]}'$, $b_{[n]}''$ and
$b_{[n]}'''$ are obtained from $b_{[n]}$ by possible changes in the
middle term, corresponding to the configuration $1^n0$, $1^n1$ and
$1^n2$, respectively. The reason is that all other configurations
contain 0 or 2 in a position before $n+1$ and therefore the changes
made by the automorphism $b$ are already accounted for in the
sequence $b_{[n]}$. Since $b(1^n0) = 1^n2$, $b(1^{n+1}0) =
1^{n+1}2$, and $b(1^n2) = 1^n0$ and the distance between $1^n0$ and
$1^n2=b(1^n0)$ along the $ca$ path is $(3^n-1)/2 + 1 + (3^n-1) + 1 +
(3^n-1)/2 = 3^{n+1}$, we see that the only change is that the middle
term in $b_{[n]}''$ is multiplied by 3. Similarly, the $3^n$ size
blocks of the $b$-change index sequence satisfy a relation of the
form $b_{[n+1]} = b_{[n]}'b_{[n]}''b_{[n]}'''$. However, the changes
in $b(1^n0) = 1^n2$ and $b(1^{n+1}0) = 1^{n+1}2$ are in the positive
direction, while the change in $b(1^n2) = 1^n0$ is in the negative
direction (we are just traveling along the same $b$ edge as in
$b(1^n0) = 1^n2$ but in the opposite direction).
\end{proof}


\section{Optimal configurations in Hanoi Towers Problem}

\begin{figure}[!ht]
\begin{center}
 \begin{tabular}{ccc}
 $
 \SelectTips{cm}{}
 \xymatrix@C=35pt{
  \A_D: &
  *++[o][F-]{{\scriptstyle()}}\ar@(ul,l)_{0,1}\ar@{}[r]|<<<{id}
  &
  *++[o][F-]{{\scriptstyle(01)}}\ar@/^1pc/[l]^{0}\ar@/_1pc/[l]_{1}\ar@{}[l]|<<{f}
  &
  *++[o][F-]{{\scriptstyle()}} \ar@{->}[l]^{1}\ar@(ur,r)^{0}\ar@{}[l]_<<{g}
 }
 $
 &&
 $
 \SelectTips{cm}{}
 \xymatrix@C=35pt{
  \A_{L_2}: &
  *++[o][F-]{{\scriptstyle()}}\ar@(ul,l)_{0} \ar@/_1pc/[r]_{1} \ar@{}[r]|<<<{\lambda_0}
  &
  *++[o][F-]{{\scriptstyle(01)}} \ar@(ur,r)^{1} \ar@/_1pc/[l]_{0} \ar@{}[l]|<<{\lambda_1}
  }
  $
  \end{tabular}
\caption{Two binary invertible transducers: $\A_D$ and
$\A_{L_2}$}\label{ad}
\end{center}
\end{figure}
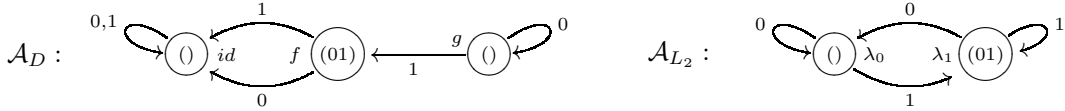

Define a matrix $M_n$ of size $2^n \times n$ with entries in $X_2$
by
\[
 M_1 = \begin{bmatrix} 0 \\ 1 \end{bmatrix}, \hspace{3cm}
 M_{n+1} = \begin{bmatrix} M_n & 0_n \\ M_n^R & 1_n \end{bmatrix},
\]
where the matrix $M_n^R$ is obtained from the matrix $M_n$ by
flipping $M_n$ along the horizontal axis, and $0_n$ and $1_n$ are
column vectors with $2^n$ entries equal to 0 and 1, respectively.
The $2^n$ rows of the matrix $M_n$ represent a binary Gray code of
length $n$. Denote the infinite limit matrix $\underset{n \to
\infty}{\lim} M_n$ by $M$. We observe that the successive rows in
$M$ are obtained from each other by applying the binary tree
automorphism $f$ at odd steps and the automorphism $g$ at even
steps, where $f$ and $g$ are given by the invertible transducer
$\A_D$ given in Figure\ref{ad}. The self-similar group $G(\A_D)$
defined by $\A_D$ and generated by $f$ and $g$ is the infinite
dihedral group $D_\infty$.

\begin{proposition}
For $j \geq 0$, define $s_{2j} = (gf)^j$ and $s_{2j+1} = f(gf)^j$.
Let $m_i$ denote the $i$-th row in the matrix $M$. Then
\[ m_i = s_i(m_0). \]
\end{proposition}

Consider the transducer in the right half of Figure~\ref{at}. It is
known~\cite{grigorchuk-z:l2} (see
also~\cite{silva-s:lamplighter,bartholdi-s:bs}) that the group
$G(\A_{L_2})$ is the lamplighter group $L_2$ which is the wreath
product of the cyclic group of order 2 (representing a switch) and
the infinite cyclic group (representing moves between consecutive
lamps). The realization of the lamplighter group $L_2$ by the
automaton $\A_{L_2}$ was used by Grigorchuk and
\.Zuk~\cite{grigorchuk-z:l2} to calculate the spectrum of the Markov
operator on the Cayley graph of $L_2$, which then lead to the
solution of Strong Atiyah Conjecture in~\cite{grigorchuk-al:atiyah}.

\begin{proposition}
For $i = 0,2^n-1$, the row $i$ word $m_i(n)$ in the matrix $M_n$
(the $i$-th Gray code word of length $n$) is equal to
$\lambda_0([i]_2^R)^R$, where $R$ denotes word reversion and $[i]_2$
is the length $n$ representative of $i$ (including leading zeros, if
necessary).
\end{proposition}

We can define a variation on the notion of transducer integer
sequences as sequences that can be obtained from transducers by
reading the input starting from the most significant digit (and
interpreting the output as starting from the most significant
digit). Call these sequences SF transducer sequences (for
significant first). Since the sequence of binary Gray code words can
be obtained by feeding the binary representations of integers, most
significant digit first, into $\A_{L_2}$ starting at $\lambda_0$, we
see that the sequence A003188 of integers
\[ 0 , 1, 3, 2, 6, 7, 5, 4, \dots \]
represented by the binary Gray code words is a SF binary transducer
sequence. On the other hand, this sequence is not a binary
transducer sequence. This is clear since the transformation $(i)_2
\mapsto m_i$ does not preserve prefixes. Namely $0^\infty \mapsto
0^\infty$, while $010^\infty \mapsto 110^\infty$.

We offer two 2 to 3 transducers each of which generates all the
configurations on the geodesic lines between the regular
configurations $0^n$, $1^n$ and $2^n$ (depending on chosen initial
state). The first one uses the order prescribed by the binary Gray
code, while the other uses the natural order.

\begin{proposition}
The transducer $\mathcal O_H$ in Figure~\ref{optimal} generates the
optimal configurations in Hanoi Towers Problem. More precisely, for
$x,y \in\{0,1,2\}$, $x \neq y$, starting at state $t_{xy}$, and
feeding the reversal $m_i(n)^R$ of the length $n$ row $i$ binary
Gray code word from $M_n$ into $\mathcal O_H$ produces the reverse
of the length $n$ ternary word representing the unique $n$ disk
configuration at distance $i$ along the geodesic from $x^n$ to $y^n$
in $\Gamma_n$.
\end{proposition}
\begin{figure}[!ht]
 \[
  \SelectTips{cm}{}
  \xymatrix@C=75pt@R=55pt{
  {\mathcal O_H:}&
  *++[o][F-]{}
  \ar@/_1pc/[r]^{1/1} \ar@/^1pc/[d]^{0/0} \ar@{}[l]|<<<<{\textstyle t_{01}}
  &
  *++[o][F-]{}
  \ar@/_1pc/[r]^{1/2} \ar@/^1pc/[d]^{0/1} \ar@{}[l]_<<<{\textstyle t_{12}}
  &
  *++[o][F-]{}
  \ar@/_2pc/[ll]_{1/0} \ar@/^1pc/[d]^{0/2} \ar@{}[r]|<<<<{\textstyle t_{20}}
  &
 \\
 &
 *++[o][F-]{}
  \ar@/_2pc/[rr]_{1/2} \ar@/^1pc/[u]^{0/0} \ar@{}[l]|<<<<{\textstyle t_{02}}
  &
  *++[o][F-]{}
  \ar@/_1pc/[l]^{1/0} \ar@/^1pc/[u]^{0/1} \ar@{}[l]^<<<{\textstyle t_{10}}
  &
  *++[o][F-]{}
  \ar@/_1pc/[l]^{1/1} \ar@/^1pc/[u]^{0/2} \ar@{}[r]|<<<<{\textstyle t_{21}}
  &
 }
 \]
 \caption{A ternary transducer generating/recognizing optimal configurations}\label{optimal}
 \end{figure}
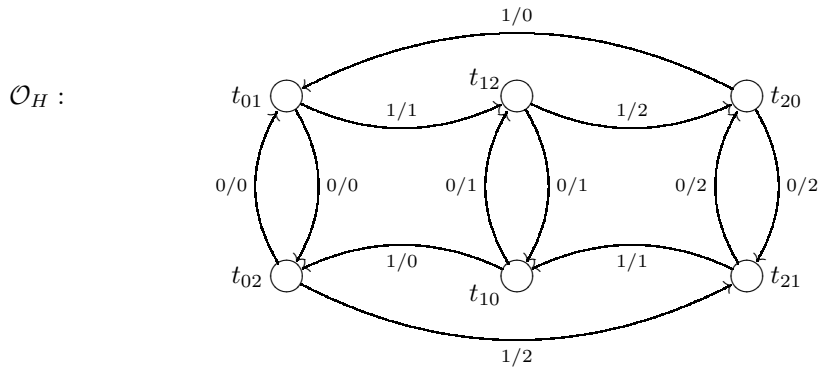
\begin{proof}
For any permutation $x,y,z$ of the three letters in $X_3$ the states
of the transducer $\mathcal O_H$ have (as tree morphisms) the
decomposition
\[
 t_{xy} =  \pi_{xy} (t_{xz},t_{yz}),
\]
where $\pi_{xy} = \left(\begin{smallmatrix}  0 & 1 \\ x & y
\end{smallmatrix}\right)$.

It is well known that the unique geodesic path of length $2^n-1$
from $x^n$ to $y^n$ connects $x^n$ to $z^{n-1}x$ in the first
$2^{n-1}-1$ steps, then in the next step the largest disk is moved
to get the configuration $z^{n-1}y$ and then the last $2^{n-1}-1$
steps are used to connect $z^{n-1}y$ to $y^n$.

Since we want to use Gray code words to describe the configurations
along the way, we observe that in the first part of the geodesic
from $x^n$ to $y^n$ (corresponding to the last digit in the Gray
code being 0) the last digit in the reached configurations is $x$,
while in the second part (corresponding to the last digit in the
Gray code being 1) the last digit in the reached configuration is
$y$. This explains the root transformations in the above
decomposition.

As for the sections, in the first part of the geodesic (last digit 0
in the Gray code) the configurations corresponding to the first
$n-1$ digits describe the path from $x^{n-1}$ to $z^{n-1}$, while in
the second part (last digit 1 in the Gray code) the configurations
corresponding to the first $n-1$ digits describe the path from
$z^{n-1}$ to $y^{n-1}$ in the natural order and from $y^{n-1}$ to
$z^{n-1}$ in the Gray code word order (because of the flip in the
second half of the Gray code). Thus the section at 0 is $t_{xz}$ and
the section at 1 is $t_{yz}$.
\end{proof}

It is apparent from the above proof that the following is also true.

\begin{figure}[!ht]
 \[
  \SelectTips{cm}{}
  \xymatrix@C=30pt@R=55pt{
  {\mathcal O_H':}
  &
  *++[o][F-]{}
  \ar@/^1pc/[rr]^{0/0} \ar@{}[ddrr]|<<<<<{\textstyle q_{01}}
  \ar@/^1pc/[dl]^{1/1}
  &&
  *++[o][F-]{}
  \ar@/^1pc/[dr]^{1/2} \ar@{}[ddll]|<<<<<{\textstyle q_{02}}
  \ar@/^1pc/[ll]^{0/0}
  &
  \\
  *++[o][F-]{}
  \ar@/^1pc/[ur]^{1/1} \ar@{}[r]|<<<<<<{\textstyle q_{21}}
  \ar@/^1pc/[dr]^{0/2}
  &&&&
  *++[o][F-]{}
  \ar@/^1pc/[dl]^{0/1} \ar@{}[l]|<<<<<{\textstyle q_{12}}
  \ar@/^1pc/[ul]^{1/2}
  \\
  &
  *++[o][F-]{}
  \ar@/^1pc/[ul]^{0/2} \ar@{}[uurr]|<<<<<{\textstyle q_{20}}
  \ar@/^1pc/[rr]^{1/0}
  &&
  *++[o][F-]{}
  \ar@/^1pc/[ll]^{1/0} \ar@{}[uull]|<<<<<{\textstyle q_{10}}
  \ar@/^1pc/[ur]^{0/1}
  &
 }
 \]
 \caption{A ternary transducer generating/recognizing optimal configurations}\label{optimal2}
 \end{figure}
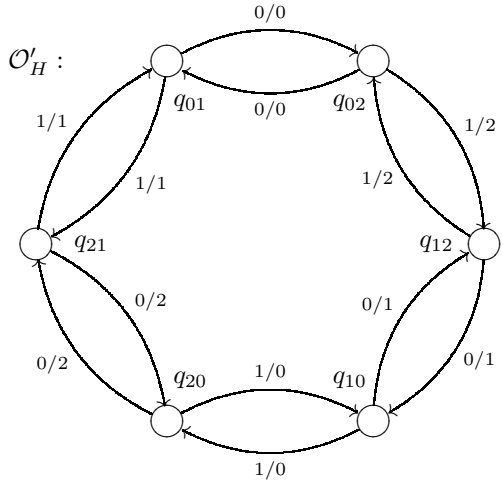

\begin{proposition}
The transducer $\mathcal O_H'$ in Figure~\ref{optimal2} generates
the optimal configurations in Hanoi Towers Problem. More precisely,
for $x,y \in\{0,1,2\}$, $x \neq y$, starting at state $q_{xy}$, and
feeding the reversal $[i]_2^R$ of the length $n$ binary
representative of $i$ (including leading 0's if needed) into
$\mathcal O_H'$ produces the reverse of the length $n$ ternary word
representing the unique $n$ disk configuration at distance $i$ along
the geodesic from $x^n$ to $y^n$ in $\Gamma_n$.
\end{proposition}

\begin{proof}
Observe that, for any permutation $x,y,z$ of the three letters in
$X_3$ the states of the transducer $\mathcal O_H'$ have (as tree
morphisms) the decomposition
\[
 q_{xy} =  \pi_{xy} (q_{xz},q_{zy}),
\]
where $\pi_{xy} = \left(\begin{smallmatrix}  0 & 1 \\ x & y
\end{smallmatrix}\right)$. This is precisely the decomposition that
corresponds to the natural order in the previous proof.
\end{proof}

The automaton $\mathcal O_H'$, started at $q_{01}$, generates the
sequence A055661
\[ 0,1,7,8,17,15,12,13,\dots, \]
but only when all input words are adjusted by leading zeros to have
odd length, and it gives the sequence
\[ 0,2,5,4,22,21,24,26,\dots, \]
which does not appear in The On-Line Encyclopedia of Integer
Sequences (as of December 2006), when the input words are adjusted
to have even length. In fact, the former sequence records the
integers whose ternary representations give the configurations in
the Hanoi Towers Problem on the geodesic line in the infinite
Schreier graph $\Gamma_{0^\infty}$ determined by applying repeatedly
the automorphisms $a$, $b$ and $c$ (in that order) and the latter
records the integers whose ternary representations give the
configurations on the geodesic line in $\Gamma_{0^\infty}$
determined by applying repeatedly the automorphisms $b$, $a$ and $c$
(in that order). There is nothing strange in this split, since it is
known that the optimal solution transferring disks to peg 1 follows
different paths depending on the parity of the number of disks.

By flipping the input and the output symbol in the automata
$\mathcal O_H$ and $\mathcal O_H'$ we obtain two automata that can
be used to recognize the configurations on the geodesic lines
between $0^\infty$, $1^\infty$ and $2^\infty$ and encode them either
by using the Gray code words or binary representations.

More generally, when $\A=(Q,X_{k_1},X_{k_2},\tau,\pi)$ is injective
transducer one can define a partial inverse transducer $\A^{-1} =
(Q^{-1}, X_{k_2},X_{k_1},\tau^{-1},\pi^{-1})$ in which
$Q=\{q^{-1}\mid q\in Q\}$, and $\tau:Q^{-1} \times X_{k_2} \to
Q^{-1}$ and $\pi:Q^{-1} \times X_{k_2} \to X_{k_1}$ are partial
maps, defined by $\tau^{-1}(q^{-1},y) = p$ and
$\pi^{-1}(p^{-1},y)=x$ whenever $\tau(q,x) = p$ and
$\pi^{-1}(p,x)=y$.

\begin{proposition}
The inverse transducer $\mathcal O_H^{-1}$, recognizes the optimal
configurations in Hanoi Towers Problem. More precisely, starting at
the inverse state $t_{xy}^{-1}$, $x,y\in X_3$, $x \neq y$, and
feeding ternary words of length $n$ into the inverse transducer
$\mathcal O_H^{-1}$, only reversals of ternary words representing
configurations on the geodesic from $x^n$ to $y^n$ in $\Gamma_n$ are
read entirely by the transducer and, for such configurations, the
reversal of the corresponding binary Gray code word of length $n$ is
produced in the output.
\end{proposition}

\begin{proposition}
The inverse transducer $\mathcal O_H'^{-1}$, recognizes the optimal
configurations in Hanoi Towers Problem. More precisely, starting at
the inverse state $q_{xy}^{-1}$, $x,y\in X_3$, $x \neq y$, and
feeding ternary words of length $n$ into the inverse transducer
$\mathcal O_H'^{-1}$, only reversals of ternary words representing
the configurations on the geodesic from $x^n$ to $y^n$ in $\Gamma_n$
are read entirely by the transducer and, for such configurations,
the output represents reversals of the binary representation of the
distance to $x^n$.
\end{proposition}

For example, the configuration $10021$ is not accepted starting from
the state $q_{01}^{-1}$ (after it is fed into $\mathcal O_H'^{-1}$
as $12001$ it stops after reading the first 4 symbols in state
$q_{20}^{-1}$ and it cannot read the last symbol). This simply means
that this configuration is not on the geodesic between $0^5$ and
$1^5$. On the other hand, $20021$ is read completely and it produces
the output $01101$, which says that the configuration $20021$ is on
the geodesic between $0^5$ and $1^5$ and its distance to $0^5$ is
$2+4+16=22$. If we read $20021$ starting at state $q_{10}^{-1}$ in
$\mathcal O_H'^{-1}$ we obtain the output $10010$, which confirms
that the configuration $20021$ is on the geodesic between $1^5$ and
$0^5$ and that its distance to $1^5$ is $1+8=9$.

\def\cprime{$'$}

\bigskip
\hrule
\bigskip

\noindent (Concerned with sequences
 A003188, A005836,
 A038500, A055661,
 A060236, A060374,
 A080846)

\end{document}